# Some Generalizations of the Chebyshev Method for Simultaneous Determination of All Roots of Polynomial Equations

## A. I. Iliev and Kh. I. Semerdzhiev

*PU "Paisii Khilendarski," ul. Tsar Asen 24, Plovdiv, 4000 Bulgaria*

**Abstract**—*Iterative methods for the simultaneous determination of all roots of an equation are dis-cussed. The multiplicities of the roots are assumed to be known in advance. The methods are proved to have a cubical rate of convergence. Numerical examples are given.*

## INTRODUCTION

During the last three–four decades, a number of classical methods for the numerical solution of equa-tions have been further advanced, which was prompted by their computer implementations and especially by using computers with parallel processors. Historically, it is interesting to note that Weierstrass predicted such methods in [1].

In the early 1960s, Newton's method was generalized in [2]; Obreshkov's method, in [3]; and the Che-byshev method, in [4] and later in [5–7]. Note that a number of generalizations of these methods to nonal-gebraic (trigonometric, exponential, and generalized) equations were obtained in [8]. They concern the case of simple roots or roots with arbitrarily given multiplicities. In [6, 7], a generalization of the Chebyshev method was developed in order to simultaneously determine all roots of the generalized polynomial $\sum_{k=0}^{n} a_k \boldsymbol{j}_k(x)$, where $a_k$ $(k = \overline{0,n})$ are given numbers and $\left\{\boldsymbol{j}_k(x)\right\}_{k=0}^{n}$ is a given Chebyshev system of basis functions. However, the Chebyshev–Makrelov method applies only to the case of simple roots. More-over, this method is rather laborious, because each iteration step involves the calculation of large-order determinants.

In this paper, we consider only algebraic, trigonometric, and exponential systems of basis functions, which are most frequently used in practice, instead of examining the general case of systems of basis func-tions. However, the resulting new methods are more effective (much less laborious), have a cubical rate of convergence, and are intended for the simultaneous determination of all roots of arbitrary multiplicities.

## 1. ALGEBRAIC POLYNOMIALS

Let an algebraic polynomial

$$A_n(x) \equiv x^n + a_1 x^{n-1} + \ldots + a_n, \tag{1.1}$$

be given. All its roots are assumed to be real. Let $x_1, \ldots, x_m$ be exact roots of (1.1) whose multiplicities are $a_1, \ldots, a_m$, respectively, $\left(a_1 + \ldots + a_m = n\right)$; and let $x_1^{[k]}, \ldots, x_m^{[k]}$ be $k$th approximations of the roots determined by the iteration

$$x_i^{[k+1]} = x_i^{[k]} - a_i \frac{A_n\left(x_i^{[k]}\right)}{A_n'\left(x_i^{[k]}\right)} \left[1 + \frac{A_n\left(x_i^{[k]}\right)}{A_n'\left(x_i^{[k]}\right)} \frac{Q_i'^{[k]}\left(x_i^{[k]}\right)}{Q_i^{[k]}\left(x_i^{[k]}\right)}\right], \quad i = \overline{1,m}, \quad k = 0,1,\ldots, \tag{1.2}$$

where

$$Q_i^{[k]}(x) = \prod_{j=1,\, j\neq i}^{m} \left(x - x_j^{[k]}\right)^{a_j}.$$

**Theorem 1.** *Let* $q$, $c$, *and* $d \overset{\text{def}}{=} \min_{i\neq j}\left|x_i - x_j\right|$ *be real numbers such that*

$0 < q < 1$, $c > 0$, $d - 2c > 0$, *and* $c^2(n - a) < (a \, d - 2nc)(d - 2c)$, *where* $i = \overline{1,m}$.

*Suppose that the initial approximations* $x_i^{[0]}$ *to the roots* $x_i$ $(i = \overline{1,m})$ *of equation* (1.1) *are*

*such that* $\left|x_i^{[0]} - x_i\right| \leq cq$, $i = \overline{1,m}$. *Then,* $\left|x_i^{[k]} - x_i\right| \leq cq^{3^k}$ $(i = \overline{1,m})$ *for any positive*

*integer* $k$.

**Proof.** We employ induction on the iteration number $k$. The conditions of the theorem imply that its assertion is true for $k = 0$. Assume that $\left|x_i^{[k]} - x_i\right| \leq cq^{3^k}$ $(i = \overline{1,m})$ $k \geq 0$. It is easy to see that

$$\frac{Q_i'^{[k]}\left(x_i^{[k]}\right)}{Q_i^{[k]}\left(x_i^{[k]}\right)} = \sum_{j=1,\, j\neq i}^{m} \frac{a_j}{x_i^{[k]} - x_j^{[k]}} \quad, i = \overline{1,m}, \tag{1.3}$$

$$\frac{A_n'\left(x_i^{[k]}\right)}{A_n\left(x_i^{[k]}\right)} = \sum_{j=1}^{m} \frac{a_j}{x_i^{[k]} - x_j} \quad, i = \overline{1,m}. \tag{1.4}$$

By using (1.3) and (1.4), method (1.2) can be transformed into

$$x_i^{[k+1]} - x_i = x_i^{[k]} - x_i - a_i\left[\frac{a_i}{x_i^{[k]} - x_i} + \sum_{j=1,\, j\neq i}^{m} \frac{a_j}{x_i^{[k]} - x_j}\right]^{-1}$$

$$- a_i \sum_{j=1,\, j\neq i}^{m} \frac{a_j}{x_i^{[k]} - x_j^{[k]}}\left[\frac{a_i}{x_i^{[k]} - x_i} + \sum_{j=1,\, j\neq i}^{m} \frac{a_j}{x_i^{[k]} - x_j}\right]^{-2},$$

or

$$x_i^{[k+1]} - x_i = x_i^{[k]} - x_i - a_i\left(x_i^{[k]} - x_i\right)/P_i^{[k]} -$$

$$- \left(x_i^{[k]} - x_i\right)^2 a_i \sum_{j=1,\, j\neq i}^{m} \frac{a_j}{x_i^{[k]} - x_j^{[k]}}\left(P_i^{[k]}\right)^{-2} \quad, i = \overline{1,m}, \tag{1.5}$$

where $P_i^{[k]}$ denotes the expression

$$a_i + \left(x_i^{[k]} - x_i\right)\sum_{j=1,\, j\neq i}^{m} \frac{a_j}{x_i^{[k]} - x_j}.$$

Extracting the factor $\left[\left(x_i^{[k]} - x_i\right)/P_i^{[k]}\right]^2$ from (1.5), we obtain

$$x_i^{[k+1]} - x_i = \left[\left(x_i^{[k]} - x_i\right)/P_i^{[k]}\right]^2\left[P_i^{[k]}\sum_{j=1,\, j\neq i}^{m} \frac{a_j}{x_i^{[k]} - x_j} - a_i\sum_{j=1,\, j\neq i}^{m} \frac{a_j}{x_i^{[k]} - x_j^{[k]}}\right],$$

or

$$x_i^{[k+1]} - x_i = \left(x_i^{[k]} - x_i\right)^2 \left[\left(x_i^{[k]} - x_i\right)\left(\sum_{j=1,\ j\neq i}^{m} \frac{\boldsymbol{a}_j}{x_i^{[k]} - x_j}\right)^2\right.$$
$$\left. + \boldsymbol{a}_i \sum_{j=1,\ j\neq i}^{m} \frac{\boldsymbol{a}_j\left(x_j - x_j^{[k]}\right)}{\left(x_i^{[k]} - x_j\right)\left(x_i^{[k]} - x_j^{[k]}\right)}\right] / \left(P_i^{[k]}\right)^2. \tag{1.6}$$

The following estimates hold:

$$\left|x_i^{[k]} - x_j\right| \geq \left|x_i - x_j\right| - \left|x_i - x_i^{[k]}\right| \geq d - cq^{3^k} > d - c > d - 2c \ ,$$

$$\left|x_i^{[k]} - x_j^{[k]}\right| \geq \left|x_i^{[k]} - x_j\right| - \left|x_j - x_j^{[k]}\right| \geq d - c - cq^{3^k} > d - 2c \quad ,i,j = \overline{1,m}\ ,i \neq j.$$

From (1.6), we find

$$\left|x_i^{[k+1]} - x_i\right| \leq \left(cq^{3^k}\right)^3 \left[\left(n - \boldsymbol{a}_i\right)^2 - \boldsymbol{a}_i\left(n - \boldsymbol{a}_i\right)\right]\left[\boldsymbol{a}_i\left(d - 2c\right) - c\left(n - \boldsymbol{a}_i\right)\right]^{-2} < cq^{3^{k+1}} \ ,i = \overline{1,m}\ ,$$

which completes the proof of Theorem 1.

**Example 1.** Consider the algebraic polynomial

$$A_6(x) = (x + 2)^2 (x - 1)(x - 3)^3$$

with the initial values $x_1^{[0]} = -3$, $x_2^{[0]} = 0.1$, and $x_3^{[0]} = 4$. Four iterations (1.2) yield the roots of the polyno-mial, which are accurate to 18 decimal digits (see Table 1).

**Table 1**

| $k$ | $x_1^{[k]}$ | $x_2^{[k]}$ | $x_3^{[k]}$ |
|---|---|---|---|
| 0 | -3.000000000000000000 | 0.100000000000000000 | 4.000000000000000000 |
| 1 | -2.074075484632669380 | 1.025215703994304140 | 3.060848242666424480 |
| 2 | -2.000104622198420050 | 0.999992663820262272 | 3.000018360022861370 |
| 3 | -2.000000000000256950 | 1.000000000000000240 | 3.000000000000001700 |
| 4 | -2.000000000000000000 | 1.000000000000000000 | 3.000000000000000000 |

## 2. TRIGONOMETRIC EQUATIONS

Consider the trigonometric polynomial

$$T_n(x) \equiv \frac{a_0}{2} + \sum_{k=1}^{n} \left[a_k \cos(kx) + b_k \sin(kx)\right], \tag{2.1}$$

where $a_n^2 + b_n^2 > 0$. Without loss of generality, we consider only the half-interval $[-\boldsymbol{p}, \boldsymbol{p}]$. Let $x_1, \ldots, x_m$ be the roots of $T_n(x)$ within this interval, whose multiplicities are $\boldsymbol{a}_1, \ldots, \boldsymbol{a}_m$, respectively, $(\boldsymbol{a}_1 + \ldots + \boldsymbol{a}_m = 2n)$; and let $x_1^{[k]}, \ldots, x_m^{[k]}$ be the $k$ th approximations to the roots determined by the following analog of method (1.2):

$$x_i^{[k+1]} = x_i^{[k]} - \boldsymbol{a}_i \frac{T_n(x_i^{[k]})}{T_n'(x_i^{[k]})} \left[1 + \frac{T_n(x_i^{[k]})}{T_n'(x_i^{[k]})} \frac{Q_i'^{[k]}(x_i^{[k]})}{Q_i^{[k]}(x_i^{[k]})}\right], \quad i = \overline{1,m}, \quad k = 0,1,\ldots, \tag{2.2}$$

where

$$Q_i^{[k]}\left(x_i^{[k]}\right) = \prod_{j=1,\ j\neq i}^{m} \sin^{\boldsymbol{a}_j}\left[\left(x_i^{[k]} - x_j^{[k]}\right)/2\right].$$

The next theorem establishes that the iterative method (2.2) has a cubical rate of convergence.

**Theorem 2.** *Let* $d \stackrel{\text{def}}{=\!=} \min_{i \neq j} \left| x_i - x_j \right|$, *and let* $c$, $q$ *and* $\boldsymbol{x}$ *be positive numbers such that* $q < 1$, $2c < \boldsymbol{x}$, $d - 2c > 0$, *and* $\max_{i \neq j} \left| x_i - x_j \right| < 2\boldsymbol{p} - 2\boldsymbol{x}$. *Suppose that* $c$ *is small enough so that*

$$c^2 \left\{ \boldsymbol{a}^2 + \frac{1}{4A^2}(2n - \boldsymbol{a})^2 + \frac{c}{4}\frac{\boldsymbol{a}}{4}(2n - \boldsymbol{a}) + \boldsymbol{a}\left[ \frac{1}{2A^2}(2n - \boldsymbol{a}) + \frac{1}{6}\frac{c}{A}(2n - \boldsymbol{a}) \right] \right\}$$

$$< \left[ \boldsymbol{a}\left( 1 - \frac{c^2}{8} \right) + \frac{c}{2}\frac{1}{A}(2n - \boldsymbol{a}) \right]^2 ,$$

*where* $A$ *denotes* $\min\left\{ \left| \sin(\boldsymbol{x}/2) \right|, \left| \sin(d/2 - c) \right| \right\}$. *Let the initial approximations* $x_i^{[0]}$ $(i = \overline{1,m})$ *be such that* $\left| x_i^{[0]} - x_i \right| \leq cq$, $i = \overline{1,m}$. *Then,* $\left| x_i^{[k]} - x_i \right| \leq cq^{3^k}$ $(i = \overline{1,m})$ *for any integer* $k \geq 0$.

The proof of this theorem follows that of Theorem 1.

For brevity, we introduce the following notation:

$$\operatorname{ctg}\frac{x_i^{[k]} - x_j^{[k]}}{2} \stackrel{\text{def}}{=\!=} t_{ij}^{[k]}, \qquad\qquad \operatorname{ctg}\frac{x_i^{[k]} - x_j}{2} \stackrel{\text{def}}{=\!=} \overline{t}_{ij}^{[k]},$$

$$\sin\frac{x_i^{[k]} - x_j^{[k]}}{2} \stackrel{\text{def}}{=\!=} s_{ij}^{[k]}, \qquad\qquad \sin\frac{x_i^{[k]} - x_j}{2} \stackrel{\text{def}}{=\!=} \overline{s}_{ij}^{[k]},$$

$$\cos\frac{x_i^{[k]} - x_j^{[k]}}{2} \stackrel{\text{def}}{=\!=} c_{ij}^{[k]}, \qquad\qquad \cos\frac{x_i^{[k]} - x_j}{2} \stackrel{\text{def}}{=\!=} \overline{c}_{ij}^{[k]}.$$

We have the presentations

$$\frac{Q_i'^{[k]}\left(x_i^{[k]}\right)}{Q_i^{[k]}\left(x_i^{[k]}\right)} = \frac{1}{2}\sum_{j=1,\,j\neq i}^{m} \boldsymbol{a}_j t_{ij}^{[k]}, \quad \frac{T_n'\left(x_i^{[k]}\right)}{T_n\left(x_i^{[k]}\right)} = \frac{1}{2}\sum_{j=1}^{m} \boldsymbol{a}_j \overline{t}_{ij}^{[k]}.$$

The iterative formula (2.2) is represented as

$$x_i^{[k+1]} - x_i = x_i^{[k]} - x_i - \boldsymbol{a}_i \frac{T_n\left(x_i^{[k]}\right)}{T_n'\left(x_i^{[k]}\right)} - \boldsymbol{a}_i \left[ \frac{T_n\left(x_i^{[k]}\right)}{T_n'\left(x_i^{[k]}\right)} \right]^2 \frac{Q_i'^{[k]}\left(x_i^{[k]}\right)}{Q_i^{[k]}\left(x_i^{[k]}\right)}$$

or, in the new notation, as

$$x_i^{[k+1]} - x_i = x_i^{[k]} - x_i - 2\boldsymbol{a}_i / \sum_{j=1}^{m} \boldsymbol{a}_j \overline{t}_{ij}^{[k]} - 2\boldsymbol{a}_i \sum_{j=1,\,j\neq i}^{m} \boldsymbol{a}_j t_{ij}^{[k]} / \left[ \sum_{j=1}^{m} \boldsymbol{a}_j \overline{t}_{ij}^{[k]} \right]^2 . \qquad (2.3)$$

We reduce the right-hand side of (2.3) to a common denominator, multiply the numerator and denominator by $\left( \overline{s}_{ij}^{[k]} \right)^2$, and extract the $i$ th term to obtain

$$x_i^{[k+1]} - x_i = \left\{ \left( x_i^{[k]} - x_i \right) \left[ \boldsymbol{a}_i \overline{c}_{ii}^{[k]} + \overline{s}_{ii}^{[k]} \sum_{j=1, j\neq i}^{m} \boldsymbol{a}_j \overline{t}_{ij}^{[k]} \right]^2 - 2\boldsymbol{a}_i \overline{s}_{ii}^{[k]} \left[ \boldsymbol{a}_i \overline{c}_{ii}^{[k]} + \overline{s}_{ii}^{[k]} \sum_{j=1, j\neq i}^{m} \boldsymbol{a}_j \overline{t}_{ij}^{[k]} \right] \right.$$

$$\left. -2\boldsymbol{a}_i \left( \overline{s}_{ii}^{[k]} \right)^2 \sum_{j=1, j\neq i}^{m} \boldsymbol{a}_j \overline{t}_{ij}^{[k]} \right\} \left[ \boldsymbol{a}_i \overline{c}_{ii}^{[k]} + \overline{s}_{ii}^{[k]} \sum_{j=1, j\neq i}^{m} \boldsymbol{a}_j \overline{t}_{ij}^{[k]} \right]^{-2}.$$

Further transformations lead to the presentation

$$x_i^{[k+1]} - x_i = \left\{ \boldsymbol{a}_i^2 \left[ \left( x_i^{[k]} - x_i \right) \left( \overline{c}_{ii}^{[k]} \right)^2 - \sin\left( x_i^{[k]} - x_i \right) \right] + \left( x_i^{[k]} - x_i \right) \left( \overline{s}_{ii}^{[k]} \right)^2 \left[ \sum_{j=1, j\neq i}^{m} \boldsymbol{a}_j \overline{t}_{ij}^{[k]} \right]^2 \right.$$

$$+ \boldsymbol{a}_i \left( x_i^{[k]} - x_i \right) \sin\left( x_i^{[k]} - x_i \right) \sum_{j=1, j\neq i}^{m} \boldsymbol{a}_j \overline{t}_{ij}^{[k]} - 2\boldsymbol{a}_i \left( \overline{s}_{ii}^{[k]} \right)^2 \sum_{j=1, j\neq i}^{m} \boldsymbol{a}_j t_{ij}^{[k]} \qquad (2.4)$$

$$\left. -2\boldsymbol{a}_i \left( \overline{s}_{ii}^{[k]} \right)^2 \sum_{j=1, j\neq i}^{m} \boldsymbol{a}_j t_{ij}^{[k]} \right\} \left[ \boldsymbol{a}_i \overline{c}_{ii}^{[k]} + \overline{s}_{ii}^{[k]} \sum_{j=1, j\neq i}^{m} \boldsymbol{a}_j \overline{t}_{ij}^{[k]} \right]^{-2}.$$

Consider the function $F(x) = x \cos^2(x/2) - \sin x$. Using the Taylor formula, we obtain $F(x) = (1/6) F'''(\boldsymbol{z}) x^3$. Then,

$$F\left( x_i^{[k]} - x_i \right) = (1/12) \left( \boldsymbol{z}_i^{[k]} \sin \boldsymbol{z}_i^{[k]} - \cos \boldsymbol{z}_i^{[k]} \right) \left( x_i^{[k]} - x_i \right)^3, \quad \boldsymbol{z}_i^{[k]} = \boldsymbol{q}_i^{[k]} \left( x_i^{[k]} - x_i \right),$$

where $0 < \boldsymbol{q}_i^{[k]} < 1$ and $i = \overline{1, m}$, and we have the estimate

$$\left| F\left( x_i^{[k]} - x_i \right) \right| \leq (1/12)(2\boldsymbol{p}+1) \left| x_i^{[k]} - x_i \right|^3 \leq \left| x_i^{[k]} - x_i \right|^3. \qquad (2.5)$$

Similarly, we have for $\Phi(x) = (x/2) \sin x - 2\sin^2(x/2)$

$$\Phi(x) = (1/24) \Phi^{(4)}(\boldsymbol{h}) x^4 : \Phi\left( x_i^{[k]} - x_i \right) = \frac{1}{24} \left( \frac{\boldsymbol{H}_i^{[k]}}{2} \sin \frac{\boldsymbol{H}_i^{[k]}}{2} - \cos \boldsymbol{H}_i^{[k]} \right) \left( x_i^{[k]} - x_i \right)^4,$$

$$\boldsymbol{h}_i^{[k]} = \boldsymbol{t}_i^{[k]} \left( x_i^{[k]} - x_i \right), \quad \boldsymbol{h}_i^{[k]} = \boldsymbol{t}_i^{[k]} \left( x_i^{[k]} - x_i \right), \quad 0 < \boldsymbol{t}_i^{[k]} < 1, \quad i = \overline{1, m},$$

$$\left| \Phi\left( x_i^{[k]} - x_i \right) \right| \leq \frac{1}{24} (\boldsymbol{p}+1) \left| x_i^{[k]} - x_i \right|^4 \leq \frac{c}{4} \left| x_i^{[k]} - x_i \right|^3. \qquad (2.6)$$

By means of $F(x)$ and $\Phi(x)$, we obtain

$$\left| x_i^{[k+1]} - x_i \right| \leq \left\{ \boldsymbol{a}_i^2 \left| F\left( x_i^{[k]} - x_i \right) \right| + \left| x_i^{[k]} - x_i \right| \left( \overline{s}_{ii}^{[k]} \right)^2 \left| \sum_{j=1, j\neq i}^{m} \boldsymbol{a}_j \overline{t}_{ij}^{[k]} \right|^2 + \boldsymbol{a}_i \left| \Phi\left( x_i^{[k]} - x_i \right) \right| \left| \sum_{j=1, j\neq i}^{m} \boldsymbol{a}_j \overline{t}_{ij}^{[k]} \right| \right.$$

$$+ \boldsymbol{a}_i \left| \frac{\left( x_i^{[k]} - x_i \right)}{2} \sin\left( x_i^{[k]} - x_i \right) \sum_{j=1, j\neq i}^{m} \boldsymbol{a}_j \overline{t}_{ij}^{[k]} - 2\left( \overline{s}_{ii}^{[k]} \right)^2 \sum_{j=1, j\neq i}^{m} \boldsymbol{a}_j t_{ij}^{[k]} \right| \right\} \left[ \boldsymbol{a}_i \left| \overline{c}_{ii}^{[k]} \right| - \left| \overline{s}_{ii}^{[k]} \right| \left| \sum_{j=1, j\neq i}^{m} \boldsymbol{a}_j \overline{t}_{ij}^{[k]} \right| \right]^{-2}. \qquad (2.7)$$

Using the equality $2\sin^2 x = 1 - \cos(2x)$ we represent $\sin x$ and $2\sin^2(x/2)$ by the Taylor formula and, thus, obtain $\sin\left( x_i^{[k]} - x_i \right) = \left( x_i^{[k]} - x_i \right) - (1/6) \cos \boldsymbol{b}_i^{[k]} \left( x_i^{[k]} - x_i \right)^3$ and $2\sin^2\left[ \left( x_i^{[k]} - x_i \right)/2 \right] = \left[ 12\left( x_i^{[k]} - x_i \right)^2 - \left( x_i^{[k]} - x_i \right)^4 \cos \boldsymbol{g}_i^{[k]} \right]/24$, where $\boldsymbol{b}_i^{[k]}$ and $\boldsymbol{g}_i^{[k]}$ are real numbers.

Hence,

$$\boldsymbol{a}\left|\frac{\left(x_i^{[k]}-x_i\right)}{2}\sin\left(x_i^{[k]}-x_i\right)\sum_{j=1,j\neq i}^m \boldsymbol{a}_j \overline{t}_{ij}^{[k]}-2\left(\overline{s}_{ii}^{[k]}\right)^2\sum_{j=1,j\neq i}^m \boldsymbol{a}_j t_{ij}^{[k]}\right|$$

$$=\boldsymbol{a}\left|\frac{\left(x_i^{[k]}-x_i\right)^2}{2}\sum_{j=1,j\neq i}^m \boldsymbol{a}_j\frac{\overline{s}_{jj}^{[k]}}{s_{ij}^{[k]}\overline{s}_{ij}^{[k]}}-\frac{1}{12}\left(x_i^{[k]}-x_i\right)^4\left[\cos \boldsymbol{b}_i^{[k]}\sum_{j=1,j\neq i}^m \boldsymbol{a}_j \overline{t}_{ij}^{[k]}-\frac{1}{2}\cos \boldsymbol{g}^{[k]}\sum_{j=1,j\neq i}^m \boldsymbol{a}_j t_{ij}^{[k]}\right]\right|.$$

Since $\left|\sin u\right|\leq|u|$, we have the estimate

$$\left|\overline{s}_{jj}^{[k]}\right|\leq\left|x_j^{[k]}-x_j\right|/2.\tag{2.8}$$

The expression in square brackets is estimated as

$$\left|\cos \boldsymbol{b}_i^{[k]}\sum_{j=1,j\neq i}^m \boldsymbol{a}_j \overline{t}_{ij}^{[k]}-\frac{1}{2}\cos \boldsymbol{g}^{[k]}\sum_{j=1,j\neq i}^m \boldsymbol{a}_j t_{ij}^{[k]}\right|\leq\left|\cos \boldsymbol{b}_i^{[k]}\right|\left|\sum_{j=1,j\neq i}^m \boldsymbol{a}_j \overline{t}_{ij}^{[k]}\right|+\frac{1}{2}\left|\cos \boldsymbol{g}^{[k]}\right|\left|\sum_{j=1,j\neq i}^m \boldsymbol{a}_j t_{ij}^{[k]}\right|$$

$$\leq\left|\sum_{j=1,j\neq i}^m \boldsymbol{a}_j \overline{t}_{ij}^{[k]}\right|+\frac{1}{2}\left|\sum_{j=1,j\neq i}^m \boldsymbol{a}_j t_{ij}^{[k]}\right|\leq\sum_{j=1,j\neq i}^m \boldsymbol{a}_j\left|\frac{\overline{c}_{ij}^{[k]}}{\overline{s}_{ij}^{[k]}}\right|+\frac{1}{2}\sum_{j=1,j\neq i}^m \boldsymbol{a}_j\left|\frac{c_{ij}^{[k]}}{s_{ij}^{[k]}}\right|\leq\sum_{j=1,j\neq i}^m \frac{\boldsymbol{a}_j}{\left|\overline{s}_{ij}^{[k]}\right|}+\frac{1}{2}\sum_{j=1,j\neq i}^m \frac{\boldsymbol{a}_j}{\left|s_{ij}^{[k]}\right|}.\tag{2.9}$$

On the other hand, we have

$$\left|x_i^{[k]}-x_j\right|\geq\left|x_i-x_j\right|-\left|x_i-x_i^{[k]}\right|\geq d-c>d-2c,$$

$$\left|x_i^{[k]}-x_j^{[k]}\right|\geq\left|x_i^{[k]}-x_j\right|-\left|x_j-x_j^{[k]}\right|\geq d-2c.$$

Moreover, we find from the assumption of Theorem 2 that $\left|x_i-x_j\right|<2\boldsymbol{p}-2\boldsymbol{x}$, $i,j=\overline{1,m}$.

Taking into account that $2c<\boldsymbol{x}$, we obtain an upper estimate:

$$\left|x_i^{[k]}-x_j\right|\leq\left|x_i^{[k]}-x_i\right|+\left|x_i-x_j\right|<2\boldsymbol{p}-2\boldsymbol{x}+cq^{3^k}<2\boldsymbol{p}-2\boldsymbol{x}+c<2\boldsymbol{p}-\boldsymbol{x},$$

$$\left|x_i^{[k]}-x_j^{[k]}\right|\leq\left|x_i^{[k]}-x_i\right|+\left|x_j^{[k]}-x_j\right|+\left|x_i-x_j\right|<2\boldsymbol{p}-2\boldsymbol{x}+2cq^{3^k}<2\boldsymbol{p}-2\boldsymbol{x}+2c<2\boldsymbol{p}-\boldsymbol{x}.$$

Thus, we have the inequalities

$$\frac{d}{2}-c<\left|\frac{x_i^{[k]}-x_j}{2}\right|<\boldsymbol{p}-\frac{\boldsymbol{x}}{2}\quad,\quad\frac{d}{2}-c<\left|\frac{x_i^{[k]}-x_j^{[k]}}{2}\right|<\boldsymbol{p}-\frac{\boldsymbol{x}}{2}$$

and, hence,

$$\left|s_{ij}^{[k]}\right|>A\quad,\quad\left|\overline{s}_{ij}^{[k]}\right|>A,\quad i,j=\overline{1,m}.\tag{2.10}$$

Then, returning to (2.8), we find that

$$\left|\cos \boldsymbol{b}_i^{[k]}\sum_{j=1,j\neq i}^m \boldsymbol{a}_j \overline{t}_{ij}^{[k]}-\frac{1}{2}\cos \boldsymbol{g}^{[k]}\sum_{j=1,j\neq i}^m \boldsymbol{a}_j t_{ij}^{[k]}\right|\leq 2\frac{\left(2n-\boldsymbol{a}\right)}{A}.\tag{2.11}$$

Let us estimate the denominator in (2.7) from below. Since $\cos x=1-\left(x^2/2\right)\cos\left(\boldsymbol{d}x\right)$ for $0<\boldsymbol{d}<1$ and $\left|\left(x_i^{[k]}-x_i\right)/2\right|<c/2$, we have $\left|\overline{c}_{ii}^{[k]}\right|>1-c^2/8$ and, hence,

$$\left[\boldsymbol{a}_i\left|\overline{c}_{ii}^{[k]}\right|-\left|\overline{s}_{ii}^{[k]}\right|\left|\sum_{j=1,j\neq i}^m \boldsymbol{a}_j t_{ij}^{[k]}\right|\right]^2>\left[\boldsymbol{a}_i\left(1-\frac{c^2}{8}\right)-\frac{c}{2}\frac{1}{A}\left(2n-\boldsymbol{a}_i\right)\right]^2.\tag{2.12}$$

Combined with (2.4)–(2.6) and (2.8)–(2.12), relation (2.7) yields

$$\left| x_i^{[k+1]} - x_i \right| \le \left| x_i^{[k]} - x_i \right|^3 \left\{ \boldsymbol{a}_i^2 + \frac{1}{4}\left( \frac{2n - \boldsymbol{a}_i}{A} \right)^2 + \boldsymbol{a}_i \frac{c}{4}\left[ \frac{1}{A}\left( 2n - \boldsymbol{a}_i \right) \right] \right.$$

$$\left. + \boldsymbol{a}_i \left[ \frac{1}{2}\left( 2n - \boldsymbol{a}_i \right)\frac{1}{A^2} + \frac{c}{6}\frac{\left( 2n - \boldsymbol{a}_i \right)}{A} \right] \right\} \left[ \boldsymbol{a}_i\left( 1 - \frac{c^2}{8} \right) + \frac{c}{2}\frac{1}{A}\left( 2n - \boldsymbol{a}_i \right) \right]^{-2}$$

$$\le c q^{3^{k+1}} c^2 \left\{ \boldsymbol{a}_i^2 + \frac{1}{4A^2}\left( 2n - \boldsymbol{a}_i \right)^2 + \frac{c}{4}\frac{\boldsymbol{a}_i}{A}\left( 2n - \boldsymbol{a}_i \right) \right.$$

$$\left. + \boldsymbol{a}_i \left[ \frac{1}{2A^2}\left( 2n - \boldsymbol{a}_i \right) + \frac{1}{6}\frac{c}{A}\left( 2n - \boldsymbol{a}_i \right) \right] \right\} \left[ \boldsymbol{a}_i\left( 1 - \frac{c^2}{8} \right) + \frac{c}{2}\frac{1}{A}\left( 2n - \boldsymbol{a}_i \right) \right]^{-2}$$

$$< c q^{3^{k+1}}, \quad i = \overline{1, m}.$$

**Example 2.** Consider the trigonometric polynomial

$$T_3(x) = \sin^3\left[ (x-1)/2 \right] \sin^2\left[ (x-2)/2 \right] \sin\left[ (x-2.5)/2 \right]$$

with the initial values $x_1^{[0]} = 0.2$, $x_2^{[0]} = 1.7$, and $x_3^{[0]} = 3$. Five iterations (2.2) produce the roots of $T_3(x)$ accurate to $18$ decimal digits (see Table 2).

**Table 2**

| $k$ | $x_1^{[k]}$ | $x_2^{[k]}$ | $x_3^{[k]}$ |
|---|---|---|---|
| 0 | 0.200000000000000000 | 1.700000000000000000 | 3.000000000000000000 |
| 1 | 1.024086327992702930 | 2.102113721613658320 | 2.719836743505084910 |
| 2 | 0.999943864177073621 | 1.994771659856962850 | 2.539910728921209960 |
| 3 | 0.999999999989823071 | 1.999997954513862020 | 2.501199355320121160 |
| 4 | 1.000000000000000000 | 1.999999999999989780 | 2.500000051660666960 |
| 5 | 1.000000000000000000 | 2.000000000000000000 | 2.500000000000000000 |

## 3. EXPONENTIAL EQUATIONS

Consider an exponential polynomial of degree $n$:

$$E_n(x) \equiv \frac{a_0}{2} + \sum_{k=1}^{n}\left( a_k \cosh x + b_k \sinh x \right). \tag{3.1}$$

Assume that at least one of the coefficients $a_n$ or $b_n$ is nonzero and all roots of the equation $E_n(x) = 0$ are real. Let $x_1, \ldots, x_m$ be the roots of the equation whose multiplicities are $\boldsymbol{a}_1, \ldots, \boldsymbol{a}_m$, respectively, $\left( \boldsymbol{a}_1 + \ldots + \boldsymbol{a}_m = 2n \right)$ and let $x_1^{[k]}, \ldots, x_m^{[k]}$ be the $k$ th approximations of the roots. To find the roots, we use the iterative formula

$$x_i^{[k+1]} = x_i^{[k]} - \boldsymbol{a}_i \frac{E_n\left( x_i^{[k]} \right)}{E_n'\left( x_i^{[k]} \right)}\left[ 1 + \frac{E_n\left( x_i^{[k]} \right)}{E_n'\left( x_i^{[k]} \right)}\frac{Q_i'\left( x_i^{[k]} \right)}{Q_i\left( x_i^{[k]} \right)} \right] \quad, i = \overline{1, m}, \quad k = 0,1,\ldots, \tag{3.2}$$

where

$$Q_i(x) = \prod_{\substack{j=1 \\ j \ne i}}^{m} \sinh^{\boldsymbol{a}_j}\left[ \left( x - x_j^{[k]} \right)/2 \right].$$

In a similar fashion to the proof of Theorem 2, one can prove that the iterative formula (3.2) has a cubical rate of convergence. Specifically, the following result holds.

**Theorem 3.** *Let* $\min\limits_{i \neq j}\left|x_i - x_j\right|$ *denote* $d$ *. Let* $q$ *and* $c$ *be real numbers such that*

$$1 > q > 0, \quad c > 0, \quad d - 2c > 0, \quad c\left|\sinh c\right| + \cosh c < 12,$$

$$\boldsymbol{a}_i^2 + nS^{-1}\left(\boldsymbol{a}_i c + S^{-3}\left|\sinh c\right|\right)\left|\sinh c\right| + 2nS^{-2}\cosh c < \boldsymbol{a}_i + S\cosh^{-1} c, \quad i = \overline{1, m}, \quad \text{where} \quad S$$

*denotes* $\sinh\left[(d - 2c)\,/\,2\right]$. *If the initial approximations* $x_i^{[0]}$ $(i = \overline{1, m})$ *are such that* $\left|x_i^{[0]} - x_i\right| \leq cq$ $(i = \overline{1, m})$, *then* $\left|x_i^{[k]} - x_i\right| \leq cq^{3^k}$ $(i = \overline{1, m})$ *for any* $k \in N$ .

**Example 3.** The iterative method (3.2) was applied to the simultaneous determination of all roots of the exponential polynomial

$$E_2(x) = \sinh^2\left[(x + 2)\,/\,2\right]\sinh^2\left[(x - 3)\,/\,2\right].$$

With the initial approximations $x_1^{[0]} = -1.5$ and $x_2^{[0]} = 3.4$, four iterations (3.2) yielded the roots of $E_2(x)$ accurate to $18$ decimal digits (see Table 3).

**Table 3**

| $k$ | $x_1^{[k]}$ | $x_2^{[k]}$ |
|---|---|---|
| 0 | -1.500000000000000000 | 3.400000000000000000 |
| 1 | -1.936759338912996590 | 3.015817214722672100 |
| 2 | -1.999910032597308230 | 3.000001221431438670 |
| 3 | -1.999999999999752340 | 3.000000000000000000 |
| 4 | -2.000000000000000000 | 3.000000000000000000 |